\documentclass{mfat}
\pagespan{74}{80}

\renewcommand{\theequation}{\mbox{\arabic{section}.\arabic{equation}}}

\begin{document}

\title[On a generalization of the three spectral inverse problem]
    {On a generalization of the three spectral inverse problem}

\author{O. P. Boyko}

\author{O. M. Martynyuk}

\author{V. M. Pivovarchik}
\address{South Ukrainian National Pedagogical University,
26 Staroportofrankivs'ka, Odesa, 65020, Ukraine}
\email{vpivovarchik@gmail.com}

\subjclass[2000]{Primary 34A55; Secondary 34B24}
\date{22/01/2015}
\dedicatory{On occasion of Yu. M. Berezansky anniversary}
\keywords{Sturm-Liouville equation, Dirichlet boundary condition, Neumann boundary condition, Marchenko equation,
Lagrange interpolation series, sine-type function, Nevanlinna function.}

\begin{abstract}
  We consider a generalization of the three spectral inverse problem,
  that is, for given spectrum of the Dirichlet-Dirichlet problem (the
  Sturm-Liouville problem with Dirichlet conditions at both ends) on
  the whole interval $[0,a]$, parts of spectra of the
  Dirichlet-Neumann and Dirichlet-Dirichlet problems on $[0,a/2]$ and
  parts of spectra of the Dirichlet-Newman and Dirichlet-Dirichlet
  problems on $[a/2,a]$, we find the potential of the Sturm-Liouville
  equation.
\end{abstract}

\maketitle

\section{Introduction}

The  theory of direct and inverse spectral  problems is based on
the base of classical results of Yu. M. Berezansky, V. A.
Marchenko, M. G. Krein, B. M. Levitan (see \cite{Be}, \cite{M},
\cite{Kr}, \cite{Le}).   The so-called half-inverse problem and
three spectral problem are branches of this theory.

In the present paper we show the relations between the three spectral
problem and the Hochstadt-Lieberman problem.  We use the same equation
(see equation (\ref{2.10})) and similar methods for recovering the
potential in these problems.

Uniqueness of solution for the three spectral problem was considered
in many publications (see \cite{GS}, \cite{D1}, \cite{D2}, \cite{D3},
\cite{FXW}, \cite{WW}, \cite{E}, \cite{EGNT}), the problem of
existence of its solution has been treated in \cite{Sa}, \cite{P} and
\cite{HM}. Let us describe the three spectral problem.

Let $\{\lambda_k\}_{-\infty, k\not=0}^{\infty}$,
$\{\nu_k^{(1)}\}_{-\infty, k\not=0}^{\infty}$ and
$\{\nu_k^{(2)}\}_{-\infty, k\not=0}^{\infty}$ be the spectra of the
problems
\begin{equation}
\label{1.1}
\left\{\begin{array}{c}
-y^{\prime\prime}+q(x)y=\lambda^2y, \\
y(0)=y(a)=0, \end{array}\right.
\end{equation}
\begin{equation}
\label{1.2}
\left\{\begin{array}{c}
-y^{\prime\prime}+q(x)y=\lambda^2y, \\
y(0)=y\left(\frac{a}{2}\right)=0, \end{array}\right.
\end{equation}
\begin{equation}
\label{1.3}
\left\{\begin{array}{c}
-y^{\prime\prime}+q(x)y=\lambda^2y, \\
y\left(\frac{a}{2}\right)=y(a)=0,  \end{array}\right.
\end{equation}
respectively, where $q\in L_2(0,a)$ is a real valued
function.

It is known \cite{P} that  $\{\nu_k^{(1)}\}_{-\infty, k\not=0}^{\infty}\cap\{\nu_k^{(2)}\}_{-\infty, k\not=0}^{\infty}=\emptyset$  if and only if $\{\lambda_k\}_{-\infty, k\not=0}^{\infty}\cap\{\nu_k^{(2)}\}_{-\infty, k\not=0}^{\infty}=\emptyset$  if and only if  $\{\lambda_k\}_{-\infty, k\not=0}^{\infty}\cap\{\nu_k^{(1)}\}_{-\infty, k\not=0}^{\infty}=\emptyset$.   In this case
this three spectra uniquely determine the potential $q$.


Let us consider also the problems
\begin{equation}
\label{1.4}
\left\{\begin{array}{c}
-y^{\prime\prime}+q(x)y=\lambda^2y, \\
y(0)=y^{\prime}\left(\frac{a}{2}\right)=0, \end{array}\right.
\end{equation}
\begin{equation}
\label{1.5}
\left\{\begin{array}{c}
-y^{\prime\prime}+q(x)y=\lambda^2y, \\
y^{\prime}\left(\frac{a}{2}\right)=y(a)=0  \end{array}\right.
\end{equation}
with the spectra denoted by  $\{\mu_k^{(1)}\}_{-\infty, k\not=0}^{\infty}$ and $\{\mu_k^{(2)}\}_{-\infty, k\not=0}^{\infty}$, respectively.

It is clear and it will be shown below that   $\{\mu_k^{(1)}\}_{-\infty, k\not=0}^{\infty}\cap\{\mu_k^{(2)}\}_{-\infty, k\not=0}^{\infty}=\emptyset$  if and only if $\{\lambda_k\}_{-\infty, k\not=0}^{\infty}\cap\{\mu_k^{(2)}\}_{-\infty, k\not=0}^{\infty}=\emptyset$  if and only if  $\{\lambda_k\}_{-\infty, k\not=0}^{\infty}\cap\{\mu_k^{(1)}\}_{-\infty, k\not=0}^{\infty}=\emptyset$.   In this case
these three spectra   $\{\lambda_k\}_{-\infty, k\not=0}^{\infty}$,  $\{\nu_k^{(1)}\}_{-\infty, k\not=0}^{\infty}$  and $\{\nu_k^{(2)}\}_{-\infty, k\not=0}^{\infty}$ uniquely determine  the potential $q(x)$ on $[0,a]$.

Let us consider another three spectral problem: given $\{\lambda_k\}_{-\infty, k\not=0}^{\infty}$,  $\{\nu_k^{(1)}\}_{-\infty, k\not=0}^{\infty}$  and $\{\mu_k^{(1)}\}_{-\infty, k\not=0}^{\infty}$, find $q$.
 Since knowledge of  $\{\mu_k^{(1)}\}_{-\infty, k\not=0}^{\infty}$  and $\{\nu_k^{(1)}\}_{-\infty, k\not=0}^{\infty}$ is equivalent to the knowledge of the projection of $q$ onto $(0,\frac{a}{2})$, this problem is nothing but the Hochstadt-Lieberman problem \cite{HL} (see also \cite{Ha}, \cite{Su}, \cite{GS}, \cite{MP},  \cite{P1}). It is known that these data uniquely determine $q$ on the whole interval $(0,a)$.

The aim of the present paper is to show that one may use  $\{\lambda_k\}_{-\infty, k\not=0}^{\infty}$  and certain parts of the spectra $\{\mu_k^{(1)}\}_{-\infty, k\not=0}^{\infty}$,   $\{\mu_k^{(2)}\}_{-\infty, k\not=0}^{\infty}$  $\{\nu_k^{(1)}\}_{-\infty, k\not=0}^{\infty}$ and  $\{\nu_k^{(2)}\}_{-\infty, k\not=0}^{\infty}$ to determine $q$. Namely,  under certain conditions, if  the spectrum $\{\lambda_k\}_{-\infty, k\not=0}^{\infty}$ is given,   the eigenvalues $\{\nu_k^{(1)}\}_{-\infty, k\not=0}^{\infty}$  and $\{\nu_k^{(2)}\}_{-\infty, k\not=0}^{\infty}$ are given not all but excluding a  finite number $2n_1$ of $\nu_k^{(1)}$s and a finite number $2n_2$ of $\nu_k^{(2)}$s then it is possible to use $2n_1$ eigenvalues $\mu_k^{(2)}$ and $2n_2$ eigenvalues $\mu_k^{(1)}$ instead to determine $q$ on $(0,a)$.

\section{Main result}
\setcounter{equation}{0} \hskip0.25in

Let us rewrite problem (\ref{1.1}) as follows

\begin{gather} \label{2.1}
-y_j^{\prime\prime}+q_j(x)y_j=\lambda^2 y_j, \quad x\in\left[0,a/2\right],
\quad  j=1,2,\\
\label{2.2} y_1(0)=0,\\
\label{2.3}
y_2(0)=0,\\
\label{2.4} y_1\left(a/2\right)=y_2\left(a/2\right),\\
\label{2.5} y_1^{\prime}(a/2)+y_2^{\prime}(a/2)=0.
\end{gather}

We rewrite problems  (\ref{1.2}) and (\ref{1.3})     as
    \begin{equation} \label{2.6}
\left\{\begin{array}{c}-y_j^{\prime\prime}+q_j(x)y_j=\lambda^2 y_j, \ \ x\in\left(0,a/2\right), \\
y_j(0)= y_j\left(a/2\right)=0, \end{array}\right. \quad j=1,2.
\end{equation}
Also we rewrite (\ref{1.4}) and (\ref{1.5}) as
    \begin{equation} \label{2.7}
\left\{\begin{array}{c}-y_j^{\prime\prime}+q_j(x)y_j=\lambda^2 y_j, \ \ x\in\left(0,a/2\right), \\
y_j(0)= y_j^{\prime}\left(a/2\right)=0, \end{array}\right. \quad
j=1,2.
\end{equation}

\medskip

\noindent
{\bf Definition 2.1.} (\cite{Y}). \ {\it ${\mathcal L}^{a}$ is the Paley-Wiener class of
entire functions of exponential type $\leq a$ which belong to
$L_2(-\infty,\infty)$ for real $\lambda$.}

\medskip

 By the Paley-Wiener theorem ${\mathcal L}^a$-functions are the Fourier images of all square summable
 functions supported on $[-a,a]$.

 Let us denote by
$s_j(\lambda,x)$ the solution of the Sturm-Liouville equation
(\ref{2.1})  which satisfies the conditions $s_j(\lambda,0)=0$,
$s_j^{\prime}(\lambda,0)=1$.
The spectra  $\{\nu_k^{(j)}\}_{-\infty, k\not=0}^{\infty}$ ($\nu_{-k}^{(j)}=-\nu_k^{(j)}$) of  problems (\ref{2.6}) coincide with the sets of zeros of the characteristic functions
$$
s_j(\lambda,\frac{a}{2})=\frac{\sin\lambda \frac{a}{2}}{\lambda}-A_j\frac{\cos\lambda
\frac{a}{2}}{\lambda^2}+ \frac{\psi_{1,j}(\lambda)}{\lambda^2},
$$
where $A_j\mathop{=}\limits^{{\operatorname{def} }}\frac{1}{2}\int_0^{\frac{a}{2}}q_j(x)\,dx$,  $\psi_{1,j}\in {\mathcal L}^{\frac{a}{2}}$.
 Moreover,
$\psi_{1,j}(0)=A_j$, $\psi_{1,j}(-\lambda)=\psi_{1,j}(\lambda)$,  otherwise $s_j(\lambda,\frac{a}{2})$ would have a pole at $\lambda=0$.

The spectra $\{\mu_k^{(j)}\}_{-\infty, k\not=0}^{\infty}$ ($\mu_{-k}^{(j)}=-\mu_k^{(j)}$) of problems (\ref{2.7}) coincide  with the set of zeros of

\begin{equation}
\label{2.9}
 s_j^{\prime}(\lambda,\frac{a}{2})=\cos\lambda \frac{a}{2}+A_j\frac{\sin\lambda
\frac{a}{2}}{\lambda}+ \frac{\psi_{2,j}(\lambda)}{\lambda},
\end{equation}
where $\psi_{2,j}\in {\mathcal L}^{\frac{a}{2}}$ and $\psi_{2,j}(0)=0$.

Let us look for the solution of problem (\ref{2.1})--(\ref{2.5}) in
the form $y_1=C_1s_1(\lambda,x)$, $y_2=C_2s_2(\lambda,x)$, where $C_j$ are
constants. Then (\ref{2.4}) and (\ref{2.5}) imply
$$ C_1s_1\left(\lambda,
a/2\right)=C_2s_2(\lambda, a/2), $$
$$ C_1s_1^{\prime}\left(\lambda, a/2\right)+C_2s_2^{\prime}(\lambda, a/2)=0.
$$
This system of equations possesses nontrivial
solution at the zeros of the characteristic function
\begin{equation}
\label{2.10}
\omega(\lambda)=s_1(\lambda,a/2)s_2^{\prime}(\lambda,a/2)+s_2(\lambda,a/2)s_1^{\prime}(\lambda,a/2).
\end{equation}
The set of zeros $\{\lambda_k\}_{-\infty, k\not=0}^{\infty}$ of this function is the spectrum of problem (\ref{2.1})--(\ref{2.5}).
Let us notice that problem (\ref{2.1})--(\ref{2.5}) is the Dirichlet-Dirichlet problem on the whole interval and therefore
$$
\omega(\lambda)=\frac{\sin\lambda a}{\lambda}-A_0\frac{\cos\lambda
a}{\lambda^2}+ \frac{\psi_{0}(\lambda)}{\lambda^2},
$$
where $A_0=A_1+A_2$,  $\psi_{0}\in {\mathcal L}^{a}$.
We introduce the following notation:
         let $\{\theta_k^2\}_{k=1}^{\infty}=:\{(\nu_k^{(1)})^2\}_{k=1}^{\infty}\cup \{(\nu_k^{(2)})^2\}_{k=1}^{\infty}$ and $\theta_{-k}=:-\theta_{k}$, $\nu_{-k}^{(j)}=-\nu_k^{(j)}$ ($j=1,2$).

The corresponding direct theorem is as follows.

\medskip

\noindent
{\bf Theorem 2.2}.  {\it Let a real-valued  $q\in L_2(0,a)$ then the spectra $\{\lambda_k\}_{-\infty, k\not=0}^{\infty}$, $\{\nu_k^{(1)}\}_{-\infty, k\not=0}^{\infty}$ and $\{\nu_k^{(2)}\}_{-\infty, k\not=0}^{\infty}$ of problems (\ref{1.1}), \  (\ref{1.2}) and  (\ref{1.3}) satisfy the conditions:

1.\ {\rm(}{\rm{see}}\cite{M}{\rm)}
\begin{equation}\label{2.13}
\begin{aligned}
\lambda_k & =\frac{\pi k}{a}+\frac{A_0}{\pi k}+\frac{\beta_k^{(0)}}{k},
\\
\nu_k^{(j)} & =\frac{2\pi k}{a}+\frac{A_j}{\pi k}+\frac{\beta_k^{(j)}}{k}, \quad
j=1,2,
\\
\mu_k^{(j)} & =\frac{\pi (2k-1)}{a}+\frac{A_j}{\pi k}+\frac{\beta_k^{(j)}}{k}, \quad j=1,2,
\end{aligned}
\end{equation}
where $A_1=\frac{1}{2}\int_0^{\frac{a}{2}}q(x)\,dx$, $A_2=\frac{1}{2}\int_{\frac{a}{2}}^aq(x)\,dx$ and $A_0=A_1+A_2$,  $\{\beta_k^{(j)}\}_{k=1}^{\infty}\in l_2$ ($j=0,1,2$).

2. All the nonzero eigenvalues $\lambda_k$, $\nu_k^{(1)}$ and $\nu_k^{(2)}$ are simple. If any of them is $0$
then it is of algebraic multiplicity 2 (\cite{M}).

3. \ {\rm(}{\rm{see}}\cite{P}{\rm)}
$$
-\infty<\lambda_1^2<\theta_1^2\leq\lambda_2^2\leq\theta_2^2\leq
\cdots
$$

4. \ {\rm(}{\rm{see}}\cite{P}{\rm)}

$$
\lambda_k^2=\theta_k^2\quad {\text {if and only if}} \quad
\lambda_k^2=\theta_{k-1}^2.
$$
}

We will use the following known  result.

\medskip

\noindent
{\bf Proposition 2.3}.  (see \cite{M}).
$$
-\infty<(\mu_1^{(j)})^2<(\nu_1^{(j)})^2<(\mu_2^{(j)})^2<(\nu_2^{(j)})^2<\cdots
$$

\medskip

\noindent
{\bf Definition 2.4}. (see e.g. \cite{PW2}).  {\it  A meromorphic function $f(z)$ is said to be an essentially positive Nevanlinna function if

(i) ${\rm Im}z {\rm Im} f(z)>0$ for all nonreal $z$,

 (ii) there exists $\beta\in\mathbb{R}$ such that $f(z)>0$ for all $z<\beta$. }

\medskip

\noindent
{\bf Proposition 2.5}. {\it  The eigenvalues $\{\lambda_k\}_{-\infty, \  k\not=0}^{\infty}$ of problem
(\ref{2.1})--(\ref{2.5})

(i) are interlaced with the union $\{\tau_k\}_{-\infty, \ k\not=0}^{\infty}\mathop{=}\limits^{{\operatorname{def}}}
\{\mu^{(1)}_k\}_{-\infty, \ k\not=0}^{\infty}\cup\{\mu^{(2)}_k\}_{-\infty, \
k\not=0}^{\infty}$
$$
-\infty<\tau_1^2\leq\lambda_1^2\leq\tau_2^2\leq\lambda_2^2\leq \cdots ,
$$

(ii)  all $\lambda_k$ are simple and for $k>0$   $\lambda_k=\tau_k$ if and only if  $\lambda_k=\tau_{k+1}$.}

\begin{proof} It is known that $\frac{s_1(\sqrt{z}, \frac{a}{2})}{s_1^{\prime}(\sqrt{z},\frac{a}{2})}$ and $\frac{s_2(\sqrt{z},\frac{a}{2})}{s_2^{\prime}(\sqrt{z},\frac{a}{2})}$ are essentially positive Nevanlinna functions.
It is known (see e.g.   \cite{PW}, Sec. 4.1) that if $f$ and $g$ are essentially positive Nevanlinna function then such is also $(f+g)$.
Therefore,
$$
\frac{\phi(\sqrt{z})}{s_1^{\prime}(\sqrt{z},\frac{a}{2})s_2^{\prime}(\sqrt{z},\frac{a}{2})}=\frac{s_1(\sqrt{z},\frac{a}{2})}{s_1^{\prime}(\sqrt{z},\frac{a}{2})}+\frac{s_2(\sqrt{z},\frac{a}{2})}{s_2^{\prime}(\sqrt{z},\frac{a}{2})}.
$$
is an essentially positive Nevanlinna function and (i) and (ii) follow by Proposition 4.3 of \cite{PW2}.
\end{proof}

\section{Inverse three spectral problem}
\setcounter{equation}{0}


In this section we will use known results on sine-type functions.

\medskip

\noindent
{\bf Definition 3.1}. {\it A function $f$ is said to be of sine-type $a$ (see \cite{LL}), if its zeros are all distinct and  there exist positive numbers $m$, $M$
and $p$ such
that
$$
me^{|{\rm Im}\lambda|a}\leq |f(\lambda)|\leq Me^{|{\rm Im}\lambda|a}
$$
for $|{\rm Im}\lambda |>p$.}

According to the statements 3 and of Theorem 2.2 each interval $(\lambda_k^2,\lambda_{k+1}^2)$ contains exactly
one element of $\{\theta_k^2\}_{k=1}^{\infty}$ and according to  Proposition 2.4  exactly one element
of $\{\tau_k^2\}_{k=1}^{\infty}$.

\medskip

\noindent
 {\bf Definition 3.2}. {\it  An interval $(\lambda_k,\lambda_{k+1})$ is said to be regular if it contains either exactly
 one element of $\{\mu_k^{(1)}\}$ and exactly one element of $\{\nu_k^{(2 )}\}$ or  it contains either exactly one
 element of $\{\mu_k^{(2)}\}$ and exactly one element of $\{\nu_k^{(1)}\}$.}

\medskip

\noindent
{\it Remark}.  It can happen that regular intervals do not exist.\

\medskip

We will use the following notation:    $$\pm\mathbb{N}=\pm 1, \pm 2, \ldots,\qquad N_j=\pm k_{1,j}, \pm k_{2,j},
\ldots , \pm k_{n_j,j}.$$

The main result of this paper is given by the following theorem.

\medskip

\noindent
{\bf Theorem 3.3}.
{\it Let  $\{\lambda_k\}_{-\infty, \  k\not=0}^{\infty}$ be the spectrum of problem (\ref{2.1})--(\ref{2.5}), $\{\nu_k^{(j)}\}_{-\infty, k\not=0}^{\infty}$ be the spectra of problems (\ref{2.6}) and let $\{\nu_k^{(1)}\}_{-\infty, k\not=0}^{\infty}\cap\{\nu_k^{(2)}\}_{-\infty, k\not=0}^{\infty}=\emptyset$.
  Let
 $\{\mu^{(2)}_k\}_{k\in  N_1}$
 ($\mu^{(2)}_{-k}=-\mu^{(2)}_k$) be eigenvalues of problem (\ref{2.7}) with $j=2$ such that $(\mu_k^{(2)})^2$  belong to the regular intervals and
 and $\{\mu^{(1)}_k\}_{k\in   N_2}$  ($\mu^{(1)}_{-k}=-\mu^{(1)}_k$) be eigenvalues of problem (\ref{2.7}) with $j=1$ such that $(\mu_k^{(1)})^2$  belong to the regular intervals.

Then  $\{\nu^{(1)}_k\}_{k\in (\pm \mathbb{N}\backslash  N_1)}$, $\{\mu^{(2)}_k\}_{k\in \pm  N_1}$,
 $\{\nu^{(2)}_k\}_{k\in (\pm \mathbb{N}\backslash  N_2)}$  and $\{\mu^{(1)}_k\}_{k\in \pm  N_2}$      uniquely determine
 the pair ($q_1$, $q_2$).
}

\begin{proof}
Using (\ref{2.13})  we find $A_j$ ($j=1,2$)
$$
A_j=\lim\limits_{k\to\infty}\pi k\left(\nu_k^{(j)}-\frac{2\pi
k}{a}\right).
$$
Let us consider the functions
$$
\omega(\lambda)=a\mathop{\prod}\limits_{k=1}^{\infty}\left(\frac{a}{\pi
k}\right)^2(\lambda_k^2-\lambda^2),
$$
\begin{equation}
\label{3.0}
\phi_1(\lambda)=\frac{a}{2}\mathop{\prod}\limits_{k\in \pm\mathbb{N}\backslash N_1}\left(\frac{a}{2\pi k}\right)^2((\nu_k^{(1)})^2-\lambda^2)\mathop{\prod}\limits_{k\in N_1}\left(\frac{a}{2\pi k}\right)^2\left((\mu_k^{(2)})^2-\lambda^2\right),
\end{equation}
\begin{equation}
\label{3.00}
\phi_2(\lambda)=\frac{a}{2}\mathop{\prod}\limits_{k\in\pm\mathbb{N}\backslash N_2}\left(\frac{a}{2\pi k}\right)^2((\nu_k^{(2)})^2-\lambda^2)\mathop{\prod}\limits_{k\in N_2}\left(\frac{a}{2\pi k}\right)^2\left((\mu_k^{(1)})^2-\lambda^2\right).
\end{equation}

 Since
 $$
s^{\prime}_j(\lambda)=\mathop{\prod}\limits_{k=1}^{\infty}\left(\frac{a}{2\pi (k-\frac{1}{2})}\right)^2
((\mu_k^{(j)})^2-\lambda^2),
$$
\begin{equation}
\label{*}
s_j(\lambda)=\frac{a}{2}\mathop{\prod}\limits_{k=1}^{\infty}\left(\frac{a}{2\pi k}\right)^2((\nu_k^{(j)})^2-\lambda^2),
\end{equation}
we conclude  keeping in mind (\ref{2.10}), (\ref{3.0}) and (\ref{3.00})
 that
the functional equation
\begin{equation}
\label{3.1}
X(\lambda)\phi_1(\lambda)+Y(\lambda)\phi_2(\lambda)=\omega(\lambda)
\end {equation}
possesses a solution
$$
\tilde{X}(\lambda)=\mathop{\prod}\limits_{k\in \pm\mathbb{N}\backslash N_1}\left(\frac{a}{\pi(2k-1)}\right)^2((\mu_k^{(2)})^2-\lambda^2)\mathop{\prod}\limits_{k\in N_1}\left(\frac{a}{\pi(2k-1)}\right)^2\left((\nu_k^{(1)})^2-\lambda^2\right),
$$
$$
\tilde{Y}(\lambda)=\mathop{\prod}\limits_{k\in \pm\mathbb{N}\backslash N_2}\left(\frac{a}{\pi(2k-1)}\right)^2((\mu_k^{(1)})^2-\lambda^2)\mathop{\prod}\limits_{k\in N_2}\left(\frac{a}{\pi(2k-1)}\right)^2\left((\nu_k^{(2)})^2-\lambda^2\right).
$$

Let us prove that equation (\ref{3.1}) possesses unique solution in the class of sine-type  functions of the form
\begin{equation}
\label{3.4}X(\lambda)=\cos\lambda \frac{a}{2}+A_2\frac{\sin\lambda
\frac{a}{2}}{\lambda}+ \frac{\tau_2(\lambda)}{\lambda},
\end{equation}
\begin{equation}
\label{3.5} Y(\lambda)=\cos\lambda \frac{a}{2}+A_1\frac{\sin\lambda
\frac{a}{2}}{\lambda}+ \frac{\tau_1(\lambda)}{\lambda},
\end{equation}
where
$\tau_j\in {\mathcal L}^{\frac{a}{2}}$.

To this end, using (\ref{2.10}) and (\ref{2.9}) we obtain
$$
s_2^{\prime}(\nu_k^{(2)},\frac{a}{2})=\frac{\omega(\nu_k^{(2)})}{s_1(\nu_k^{(2)},\frac{a}{2})}=\cos\nu_k^{(2)}\frac{a}{2}+A_2\frac{\sin\nu_k^{(2)}\frac{a}{2}}{\nu_k^{(2)}}+\frac{\psi_{2,j}(\nu_k^{(2)})}{\nu_k^{(2)}},
$$
where due to (2.13) $\{\psi_{2,j}(\nu_k^{(2)})\}_{-\infty, k\not=0}^{\infty}\in l_2$ by Lemma 1.4.3 of \cite{M} .

We look for the solution of (\ref{3.1}) in the form (\ref{3.4}), (\ref{3.5}).  Let us consider (\ref{3.1}) at the zeros of $\phi_2$ which we will use as the nodes of interpolation.
Comparing (\ref{3.0}) and (\ref{*}) we obtain
\begin{equation}
\label{3.6}
\phi_1(\lambda)=s_1(\lambda,\frac{a}{2})\mathop{\prod}\limits_{k\in N_1}\left((\mu_k^{(2)})^2
-\lambda^2\right)((\nu_k^{(1)})^2-\lambda^2)^{-1}.
\end{equation}
Thus, we conclude that      for $k\in\pm{\mathbb N}\backslash N_2 $
$$
X(\nu_k^{(2)})=\frac{\omega(\nu_k^{(2)})}{\phi_1(\nu_k^{(2)})}=\frac{\omega(\nu_k^{(2)})}{s_1(\nu_k^{(2)},\frac{a}{2})}
\mathop{\prod}\limits_{p\in N_1}\left((\mu_p^{(2)})^2-(\nu_k^{(2)})^2\right)^{-1}((\nu_p^{(1)})^2-(\nu_k^{(2)})^2)
$$
and for $k\in N_2$
$$
X(\mu_k^{(1)})=\frac{\omega(\mu_k^{(1)})}{\phi_1(\mu_k^{(1)})}=\frac{\omega(\mu_k^{(1)})}{s_1(\mu_k^{(1)},\frac{a}{2})}
\mathop{\prod}\limits_{p\in
N_1}\left((\mu_p^{(2)})^2-(\mu_k^{(1)})^2\right)^{-1}((\nu_p^{(1)})^2-(\mu_k^{(1)})^2).
$$
Hence, for $k$ large
$$
\begin{aligned}
X(\nu_k^{(2)}) & =\left(\cos\nu_k^{(2)}\frac{a}{2}+A_2\frac{\sin\nu_k^{(2)}\frac{a}{2}}{\nu_k^{(2)}}+
\frac{\psi_{2,j}(\nu_k^{(2)})}{\nu_k^{(2)}}\right)(1+O(k^{-2})
\\
& = \cos\nu_k^{(2)}\frac{a}{2}+A_2\frac{\sin\nu_k^{(2)}\frac{a}{2}}{\nu_k^{(2)}}+
\left(\cos\nu_k^{(2)}\frac{a}{2}+A_2\frac{\sin\nu_k^{(2)}\frac{a}{2}}{\nu_k^{(2)}}\right)
O(k^{-2})
\\
& +\frac{\psi_{2,j}(\nu_k^{(2)})}{\nu_k^{(2)}}(1+O(k^{-2})).
\end{aligned}
$$
Than the sequence $\{\tau_k^{(2)}\}_{-\infty, k\not=0}^{\infty}$ where
$$
\tau_k^{(2)}=: \nu_k^{(2)}\left( X(\nu_k^{(2)})-\cos\nu_k^{(2)}\frac{a}{2}-A_2\frac{\sin\nu_k^{(2)}\frac{a}{2}}{\nu_k^{(2)}}\right), \ \  k\in\pm\mathbb{N}\backslash N_2
$$
and
$$
\tau_k^{(2)}=: \mu_k^{(1)}\left( X(\mu_k^{(1)})-\cos\mu_k^{(1)}\frac{a}{2}-A_2\frac{\sin\mu_k^{(1)}\frac{a}{2}}{\mu_k^{(1)}}\right), \ \  k\in  N_2
$$
belongs to $ l_2$.

Let us assume that $\phi_2(0)\ne 0$, otherwise we can shift the spectral parameter $\lambda^2\rightarrow \lambda^2+c$. Since $\tilde{\phi_2}(\lambda)=\lambda \phi_2(\lambda) $ is a sine type function with all zeros including  $\tilde{\phi}_2(0)=0$ simple we can use the set $\{\nu_k^{(2)}\}_{k\in\pm{\mathcal N}\backslash N_2}\cup\{\mu_k^{(1)}\}_{k\in N_2}\cup\{0\}$
    of zeros of $\tilde{\phi}_2$ as the nodes of interpolation. The values of the function $\tau_2(\lambda)$ we interpolate are $\tau^{(2)}_k$ and $\tau^{(2)}(0)=0$.
Since $\{\tau_k^{(2)}\}_{-\infty, k\not=0}^{\infty}\in l_2$, according to \cite{LL}, \cite{Y}   the series
\begin{equation}
\label{3.8}
\tau_2(\lambda)\!=:\!\tilde{\phi}_2(\lambda)\mathop{\sum}\limits_{k\in\pm{\mathcal N}\backslash N_2}
\frac{\tau_k^{(2)}}{\frac{d\tilde{\phi}_2(\lambda)}{d\lambda}|_{\lambda=\nu_k^{(2)}}(\lambda-\nu_k^{(2)})}\!+\!\tilde{\phi}_2(\lambda)\mathop{\sum}\limits_{k\in N_2}\frac{\tau_k^{(2)}}{\frac{d\tilde{\phi}_2(\lambda)}{d\lambda}|_{\lambda=\mu_k^{(1)}}(\lambda-\mu_k^{(1)})}
\end{equation}
converges uniformly on any compact of $\mathbb{C}$ and in the norm of $L_2(-\infty, \infty)$ to a function in ${\mathcal L}^{a/2}$ . Since (\ref{3.8}) establishes one-to-one correspondence between $l_2$ and ${\mathcal L}^{\frac{a}{2}}$ we conclude that
$$
X(\lambda)=:\cos\lambda\frac{a}{2}+A_2\frac {\sin\lambda\frac{a}{2}}{\lambda}+\frac{\tau_2(\lambda)}{\lambda}=\tilde{X}(\lambda).
$$
It is easy to prove in the same way that $\tilde{Y}(\lambda)$ is also uniquely determined by  (\ref{3.1}) . We  identify with $\nu_k^{(1)}$ the  zero of $\tilde{X}$  which lies in the regular interval containing $\mu_k^{(2)}$ ($k\in N_1$) and  with $\nu_k^{(2)}$ the  zero of $\tilde{Y}$  which lies in the regular interval containing $\mu_k^{(1)}$ ($k\in N_2$).    It is clear that
$\{\nu_k^{(j)}\}_{-\infty, k\not=0}^{\infty}$ and $\{\mu_k^{(j)}\}_{-\infty, k\not=0}^{\infty}$ uniquely
determine $q_j(x)$. The method of recovering $q_j(x)$ by the two spectra $\{\nu_k^{(j)}\}_{-\infty, k\not=0}^{\infty}$
and $\{\mu_k^{(j)}\}_{-\infty, k\not=0}^{\infty}$ can be found in~\cite{M}.
The theorem is proved.
\end{proof}



\end{document}